\documentclass[twoside]{article}
\usepackage[english]{babel}
\usepackage[cp1251]{inputenc}
\usepackage{amssymb,amsmath,amsthm}
\begin{document}
\pagestyle{myheadings}

\markboth{On isomorphism classes and invariants...}{I.S. Rakhimov,
\ \ S.K. Said Husain}

\begin{center}
\Large {On isomorphism classes and invariants of low dimensional
complex filiform Leibniz algebras (part 1)}
\end{center}

\begin{center}
I.S. Rakhimov$^1$, \ \ S.K. Said Husain$^2$
\end{center}

$^1$isamiddin@science.upm.edu.my $\&$ risamiddin@mail.ru

$^2$skartini@science.upm.edu.my

\begin{center}
Institute for Mathematical Research (INSPEM) $\&$ Department of
Mathematics, FS,\\ 43400 UPM, Serdang, Selangor Darul Ehsan,
(Malaysia)
\end{center}

\begin{abstract}
The paper aims to investigate the classification problem of low
dimensional complex none Lie filiform Leibniz algebras. There are
two sources to get classification of filiform Leibniz algebras.
The first of them is the naturally graded none Lie filiform
Leibniz algebras and the another one is the naturally graded
filiform Lie algebras \cite{GO}. Here we do consider Leibniz
algebras appearing from the naturally graded none Lie filiform
Leibniz algebras. According to the theorem presented in \cite{AO}
this class can be splited into two subclasses. However,
isomorphisms within each class there were not investigated. In
\cite{BR} U.D.Bekbaev and I.S.Rakhimov suggested an approach to
the isomorphism problem in terms of invariants. This paper
presents an implementation of the results of \cite{BR} in low
dimensional cases. Here we give the complete classification of
complex none Lie filiform Leibniz algebras in dimensions at most 8
from the first class of the above mentioned result of \cite{AO}
and give a hypothetic formula for the number of isomorphism
classes in finite dimensional case.
\end{abstract}

\footnote{The research is supported by Grant 04-01-06 SF01-22
MOSTI (Malaysia)}



\thispagestyle{empty}
\section{Preliminaries}

{\bf Definition 1.1} An algebra $L$ over a field $K$ is called a
{\it Leibniz algebra} if it satisfies the following Leibniz
identity: $$x(yz)=(xy)z-(xz)y.$$

Let $Leib_{n}(K)$ be a subvariety of $Alg_{n}(K)$ consisting of
all $n$-dimensional Leibniz algebras over $K$. It is invariant
under the isomorphic action of $GL_n(K)$ ("transport of
structure"). Two algebras are isomorphic if and only if they
belong to the same orbit under this action. As a subset of
$Alg_{n}(K)$ the set $Leib_{n}(K)$ is specified by the system of
equations with respect to structural constants $\gamma_{ij}^{k}$:
$$\sum\limits_{\emph{l}=1}^{\emph{n}}{(\gamma_{\emph{jk}}^{\emph{l}}
\gamma_{\emph{il}}^{\emph{m}}-\gamma_{\emph{ij}}^{\emph{l}}\gamma_{\emph{lk}}^{\emph{m}}+
\gamma_{\emph{ik}}^{\emph{l}}\gamma_{\emph{lj}}^{\emph{m}})}=0.$$

It is easy to see that if the multiplication in Leibniz algebra
happens to be anticommutative then it is a Lie algebra. So Leibniz
algebras are "noncommutative" generalization of Lie algebras. As
for Lie algebras case they are well known and several
classifications of low dimensional cases have been given. But
unless simple Lie algebras the classification problem of all Lie
algebras in common remains a big problem. Yu.I.Malcev \cite{Mal}
reduced the classification of solvable Lie algebras to the
classification of nilpotent Lie algebras. Apparently the first
non-trivial classification of some classes of low-dimensional
nilpotent Lie algebra are due to Umlauf. In his thesis \cite{Um}
he presented the redundant list of nilpotent Lie algebras of
dimension less seven. He gave also the list of nilpotent Lie
algebras of dimension less than ten admitting a so-called adapted
basis (now, the nilpotent Lie algebras with this property are
called \emph{filiform Lie algebras}). It was shown by M.Vergne
\cite{Vr} the importance of filiform Lie algebras in the study of
variety of nilpotent Lie algebras laws.

Further if it is not asserted additionally all algebras assumed to
be over the field of complex numbers.

Let $L$ be a Leibniz algebra. We put: $ L^1=L,\quad
L^{k+1}=[L^k,L],\enskip k\in N.$

{\bf Definition 1.2} A Leibniz algebra $L$ is said to be {\it
nilpotent} if there exists an integer $s\in N,$ such that $L^{1}
\supset L^{2} \supset ... \supset L^{s}=\{0\}.$ The smallest
integer $s$ for which $L^{s}=0$ is called {\it the nilindex} of
$L$.

{\bf Definition 1.3} An $n$-dimensional Leibniz algebra $L$ is
said to be {\it filiform} if $dim L^i =n-i,$ where $2\le i\le n.$

{\bf Theorem 1.4 \cite{AO}} Arbitrary complex non-Lie filiform
Leibniz algebra of dimension $n+1$ obtained from naturally graded
none Lie filiform Leibniz algebra is isomorphic to one of the
following filiform Leibniz algebras with the table non zero
multiplications for basis vectors $\{e_0,e_1, \ldots, e_n\}$:

$ \mbox{ a) (The first class):}  \left\{
\begin{array}{lll}
e_{0}e_{0}=e_{2},\\

e_ie_{0}=e_i+1, \qquad \qquad \qquad \qquad \qquad
\qquad \qquad  \ \ 1\leq i\leq n-1 \\

e_{0}e_{1}]= \alpha_{3} e_{3}+ \alpha_{4} e_{4}+...+\alpha_{n-1}
e_{n-1}+
\theta e_{n}, \\

e_{j}e_{1}=\alpha_{3}e_{j+2}+
\alpha_{4}e_{j+3}+...+\alpha_{n+1-j}e_{n}, \qquad 1\leq j\leq n-2
\end{array}
\right.
$

$ \mbox{b) (The second class):}  \left\{
\begin{array}{lll}
\emph{e}_{0}\emph{e}_{0}=\emph{e}_{2},\\

e_{i}e_{0}=e_{i+1}, \qquad \qquad \qquad \qquad \qquad
\qquad \qquad  \ \ 2\leq i\leq n-1 \\

e_{0}e_{1}= \beta_{3} e_{3}+ \beta_{4}
e_{4}+...+\beta_{n} e_{n}, \\

e_{1}e_{1}=\gamma e_{n},\\

e_{j}e_{1}=\beta_{3}e_{j+2}+
\beta_{4}e_{j+3}+...+\beta_{n+1-j}e_{n}, \qquad 2\leq j\leq n-2

\end{array}
\right.$

Note that the algebras from the first class and the second class
never are isomorphic to each other.

In this paper we will deal with the first class of algebras of the
above Theorem, they will be denoted
$L(\alpha_{3},\alpha_{4},...,\alpha_{\emph{n}},\theta)$, meaning
that they are defined by parameters
$\alpha_{3},\alpha_{4},...,\alpha_{\emph{n}},\theta$. The class
here will be denoted as $FLeib_{n+1}$. As for the second class it
will be considered somewhere else.

Using the method of simplification of the basis transformations in
\cite{GO} the following criterion on isomorphism of two
$(n+1)$-dimensional filiform Leibniz algebras was given. Namely:
let $n\geq 3$.

{\bf Theorem 1.5} Two algebras $L(\alpha)$ and $L(\alpha')$ from
$FLeib_{n+1}$, where
$\alpha=(\alpha_{3},\alpha_{4},...,\alpha_{\emph{n}},\theta)$ and
$\alpha'=(\alpha'_{3},\alpha'_{4},...,\alpha'_{\emph{n}},\theta')$,
are isomorphic if and only if there exist complex numbers
$\emph{A},\emph{B}$ such that $\emph{A}(\emph{A}+\emph{B})\neq$ 0
and the following conditions hold:

$ \left\{
\begin{array}{lll}
\alpha'_{3}=\frac{A+B}{A^{2}}\alpha_{3},\\
\alpha'_{t}=\frac{1}{A^{t-1}}((A+B)\alpha_{t}- \sum
\limits_{k=3}^{t-1}(C_{k-1}^{k-2}A^{k-2}B \alpha_{t+2-k}+
C_{k-1}^{k-3}A^{k-3}B^{2} \sum \limits_{i_{1}=k+2}^{t}
\alpha_{t+3-i_{1}}\cdot\alpha_{i_{1}+1-k}+\\
C_{k-1}^{k-4}A^{k-4}B^{3}\sum\limits_{i_{2}=k+3}^{t}
\sum\limits_{i_{1}=k+3}^{i_{2}}\alpha_{t+3-i_{2}}\cdot\alpha_{i_{2}+3-i_{1}}\cdot
\alpha_{i_{1}-k}+...+\\
C_{k-1}^{1}AB^{k-2}\sum\limits_{i_{k-3}=2k-2}^{t}
\sum\limits_{i_{k-4}=2k-2}^{i_{k-3}}...
\sum\limits_{i_{1}=2k-2}^{i_{2}}\alpha_{t+3-i_{k-3}}\cdot\alpha_{i_{k-3}+3-i_{k-4}}
\cdot...\cdot\alpha_{i_{2}+3-i_{1}}\cdot\alpha_{i_{1}+5-2k}\\+B^{k-1}
\sum\limits_{i_{k-2}=2k-1}^{t}
\sum\limits_{i_{k-3}=2k-1}^{i_{k-2}}...
\sum\limits_{i_{1}=2k-1}^{i_{2}}\alpha_{t+3-i_{k-2}}\cdot\alpha_{i_{k-2}+3-i_{k-3}}
\cdot...\cdot\alpha_{i_{2}+3-i_{1}}\alpha_{i_{1}+4-2k})\cdot\alpha'_{k}), \\

\mbox{where} \ \ 4 \leq t \leq n.\\

\theta'=\frac{1}{A^{n-1}}(A\theta+B\alpha_{n}- \sum
\limits_{k=3}^{n-1}(C_{k-1}^{k-2}A^{k-2}B \alpha_{n+2-k}+
C_{k-1}^{k-3}A^{k-3}B^{2} \sum \limits_{i_{1}=k+2}^{n}
\alpha_{n+3-i_{1}}\cdot\alpha_{i_{1}+1-k}+\\
C_{k-1}^{k-4}A^{k-4}B^{3}\sum\limits_{i_{2}=k+3}^{n}
\sum\limits_{i_{1}=k+3}^{i_{2}}\alpha_{n+3-i_{2}}\cdot\alpha_{i_{2}+3-i_{1}}
\cdot\alpha_{i_{1}-k}+...+\\
C_{k-1}^{1}AB^{k-2}\sum\limits_{i_{k-3}=2k-2}^{n}
\sum\limits_{i_{k-4}=2k-2}^{i_{k-3}}...
\sum\limits_{i_{1}=2k-2}^{i_{2}}\alpha_{n+3-i_{k-3}}\cdot\alpha_{i_{k-3}+3-i_{k-4}}
\cdot...\cdot\alpha_{i_{2}+3-i_{1}}\cdot\alpha_{i_{1}+5-2k}\\+B^{k-1}
\sum\limits_{i_{k-2}=2k-1}^{n}
\sum\limits_{i_{k-3}=2k-1}^{i_{k-2}}...
\sum\limits_{i_{1}=2k-1}^{i_{2}}\alpha_{n+3-i_{k-2}}\cdot\alpha_{i_{k-2}+3-i_{k-3}}
\cdot...\cdot\alpha_{i_{2}+3-i_{1}}\cdot\alpha_{i_{1}+4-2k})\cdot\alpha'_{k}),

\end{array}
\right.
$
\\

It is not difficult to notice that the expressions for
$\alpha'_{t}$, $\theta'$ in Theorem 1.5 can be represented in the
following form:
$\alpha'_{t}=\frac{1}{A^{t-2}}\varphi_{t}(\frac{B}{A};\alpha),$
where $\alpha=(\alpha_{3},\alpha_{4},...,\alpha_{n}, \theta)$ and
$\varphi_{t}(y;z)=\varphi_{t}(y;z_{3},z_{4},...,z_{n},z_{n+1})=$
$$\begin{array}{lll} ((1+y)z_{t}- \sum \limits_{k=3}^{t-1}(C_{k-1}^{k-2}y z_{t+2-k}+
C_{k-1}^{k-3}y^{2} \sum \limits_{i_{1}=k+2}^{t} z_{t+3-i_{1}}
\cdot z_{i_{1}+1-k}+\\
\\
C_{k-1}^{k-4}y^{3}\sum\limits_{i_{2}=k+3}^{t}
\sum\limits_{i_{1}=k+3}^{i_{2}} z_{t+3-i_{2}} \cdot
z_{i_{2}+3-i_{1}} \cdot z_{i_{1}-k}+...+\\
\\
C_{k-1}^{1}y^{k-2}\sum\limits_{i_{k-3}=2k-2}^{t}
\sum\limits_{i_{k-4}=2k-2}^{i_{k-3}}...
\sum\limits_{i_{1}=2k-2}^{i_{2}} z_{t+3-i_{k-3}} \cdot
z_{i_{k-3}+3-i_{k-4}} \cdot...\cdot z_{i_{2}+3-i_{1}} \cdot
z_{i_{1}+5-2k}+\\
\\y^{k-1} \sum\limits_{i_{k-2}=2k-1}^{t}
\sum\limits_{i_{k-3}=2k-1}^{i_{k-2}}...
\sum\limits_{i_{1}=2k-1}^{i_{2}} z_{t+3-i_{k-2}} \cdot
z_{i_{k-2}+3-i_{k-3}} \cdot...\cdot z_{i_{2}+3-i_{1}}
z_{i_{1}+4-2k}) \cdot \varphi_{k}(y;z)),\\
\\
\mbox{for} \ 3 \leq t \leq n.\end{array}$$
$\theta'=\frac{1}{A^{n-2}}\varphi_{n+1}(\frac{B}{A};\alpha),
\mbox{where}\
\varphi_{n+1}(y;z)=\varphi_{n+1}(y;z_{3},z_{4},...,z_{n},z_{n+1})=$
$$\begin{array}{lll}
(z_{n+1}+y z_{n}- (1+y)\sum \limits_{k=3}^{n-1}(C_{k-1}^{k-2}y
z_{n+2-k}+ C_{k-1}^{k-3}y^{2} \sum \limits_{i_{1}=k+2}^{n}
z_{n+3-i_{1}} \cdot z_{i_{1}+1-k}+\\
\\
C_{k-1}^{k-4}y^{3}\sum\limits_{i_{2}=k+3}^{n}
\sum\limits_{i_{1}=k+3}^{i_{2}} z_{n+3-i_{2}} \cdot
z_{i_{2}+3-i_{1}}
\cdot z_{i_{1}-k}+...+\\
\\
C_{k-1}^{1}y^{k-2}\sum\limits_{i_{k-3}=2k-2}^{n}
\sum\limits_{i_{k-4}=2k-2}^{i_{k-3}}...
\sum\limits_{i_{1}=2k-2}^{i_{2}}z_{n+3-i_{k-3}} \cdot
z_{i_{k-3}+3-i_{k-4}} \cdot...\cdot z_{i_{2}+3-i_{1}}\cdot
z_{i_{1}+5-2k}\\
\\
+y^{k-1} \sum\limits_{i_{k-2}=2k-1}^{n}
\sum\limits_{i_{k-3}=2k-1}^{i_{k-2}}...
\sum\limits_{i_{1}=2k-1}^{i_{2}}z_{n+3-i_{k-2}}\cdot
z_{i_{k-2}+3-i_{k-3}} \cdot...\cdot z_{i_{2}+3-i_{1}}\cdot
z_{i_{1}+4-2k})\cdot\varphi_{k}(y;z)).\end{array}$$

For transition from the $(n+1)$-dimensional filiform Leibniz
algebra $L(\alpha)$ to the $(n+1)$-dimensional filiform Leibniz
algebra $L(\alpha')$ we will write
$\alpha'=\rho(\frac{1}{A},\frac{B}{A};\alpha)$, \ \ where
$\alpha=(\alpha_{3},\alpha_{4},...,\alpha_{n},\theta)$,

$$\rho(\frac{1}{A},\frac{B}{A};\alpha)=(\rho_{1}(\frac{1}{A},\frac{B}{A};\alpha),
\rho_{2}(\frac{1}{A},\frac{B}{A};\alpha),...,\rho_{n-1}(\frac{1}{A},\frac{B}{A};
\alpha)),$$

$$\rho_{t}(x,y;z)=x^{t}\varphi_{t+2}(y;z)\ \mbox{for} \ 1 \leq t \leq
n-2$$ \ and \

$$\rho_{n-1}(x,y;z)=x^{n-2}\varphi_{n+1}(y,z)$$

Here are the main properties of the operator $\rho$, derived from
the fact that $\rho(\frac{1}{A},\frac{B}{A};\cdot)$ is an action
of a group.

$$\begin{array}{lll}{1^{0}. \ \ \rho(1,0;\cdot) \ \  \mbox{is the identity
operator}.}\\
\\ 2^{0}. \
\rho(\frac{1}{A_{2}},\frac{B_{2}}{A_{2}};\rho(\frac{1}{A_{1}},
\frac{B_{1}}{A_{1}};\alpha))=\rho(\frac{1}{A_{1}A_{2}},
\frac{A_{1}B_{2}+A_{2}B_{1}+B_{1}B_{2}}{A_{1}A_{2}};\alpha).\\
\\ 3^{0}. \ \ \mbox{If} \ \
\alpha'=\rho(\frac{1}{A},\frac{B}{A};\alpha) \ \ \mbox{then} \ \
\alpha=\rho(A,-\frac{B}{A+B};\alpha').\end{array}$$

From here on $n$ is a positive integer. We assume that $n\geq4$
since there are complete classifications of complex nilpotent
Leibniz algebras of dimension at most four \cite{AOR}.

We first present the result of \cite{BR} that underlies our
classification result.

We consider the following presentation of $FLeib_{n+1}:$
$FLeib_{n+1}=U_{1} \cup F$, where $U_1= \{L(\alpha):
\alpha_{3}(\alpha_{4}+2\alpha_{3}^{2}) \neq 0\}$, $F= \{L(\alpha):
\alpha_{3}(\alpha_{4}+2\alpha_{3}^{2}) = 0\}.$

{\bf Theorem 1.6} \cite{BR} $i)$ Two algebras $L(\alpha)$ and
$L(\alpha')$ from $\emph{U}_1$ are isomorphic if and only if
$$
\rho_{i}(\frac{2\alpha_{3}}{\alpha_{4}+ 2
\alpha_{3}^{2}},\frac{\alpha_{4}}{2 \alpha_{3}^{2}};\alpha) =
\rho_{i}(\frac{2\alpha_{3}'}{\alpha_{4}'+ 2
\alpha_{3}'^{2}},\frac{\alpha_{4}'}{2 \alpha_{3}'^{2}};\alpha')
$$
whenever \ \ $\emph{i}=\overline{3,n-1}.$

$ii)$ For any $(a_3,a_4,...,a_{n-1})\in C^{n-3}$ there is an
algebra $L(\alpha)$ from $\emph{U}_1$ such that
$$\rho_{i}(\frac{2\alpha_{3}}{\alpha_{4}+ 2
\alpha_{3}^{2}},\frac{\alpha_{4}}{2 \alpha_{3}^{2}};\alpha)= a_i \
\ \mbox{for all} \ \ \emph{i}=\overline{3,n-1}.$$

Note that the above theorem describes the field of invariant
rational functions on $Leib_{n+1}$ under the action of the adapted
subgroup of $GL_{n+1}(\mathbf{C})$ and the part $ii)$ means the
algebraic independence of the generators.

The procedure that we are applying works by the following way:
first we present $FLeib_{n+1}$ as a disjoint union of subsets then
formulate for each subset the analogue of the above theorem.

\section{The list of algebras}

In this section we will present the list of none Lie complex
filiform Leibniz algebras from $FLeib_{n+1}$ for $n=4,5,6,7.$
Later on $\Delta_4=\alpha_4+2\alpha_3^2, \ \
\Delta_5=\alpha_5-5\alpha_3^3, \ \ \Delta_6=\alpha_6+14\alpha_3^4,
\ \ \Delta_7=\alpha_7-42\alpha_3^5, \ \
 \ \ \Theta_i=\theta-\alpha_i, \ \
i=4,5,6,7$ and the same letters $\Delta$ and $\Theta$ with
$\prime$ will denote the same expression depending on parameters
$\alpha_3', \alpha_4', \alpha_5', \alpha_6',\alpha_7', \theta'.$
Notice that $\Delta_i=\alpha_i$ $(i=4,5,6,7)$ when $\alpha_3=0.$

Let $N_n$ denote the number of isomorphism classes in dimension
$n$ (each parametric family here will be considered as a one
class).

\textbf{2.1 Dimension 5}

The class $FLeib_{5}$ can be represented as a disjoint union of
the following subsets:

$\qquad FLeib_{5}=U_{1}\bigcup U_{2}\bigcup U_{3}\bigcup
U_{4}\bigcup U_{5}\bigcup U_{6}\bigcup U_{7},$ where

$\qquad U_1=\{L(\alpha)\in FLeib_5:\alpha_3\neq 0, \Delta_4 \neq 0
\},$

$\qquad U_2=\{L(\alpha)\in FLeib_5:\alpha_3\neq 0, \Delta_4=
0,\Theta_4\neq0\} ,$

$\qquad U_{3}=\{L(\alpha )\in FLeib_{5}:\alpha _{3}\neq
0,\Delta_4=0,\Theta_4=0 \},$

$\qquad U_{4}=\{L(\alpha )\in FLeib_{5}:\alpha _{3}=0,\Delta
_{4}\neq 0,\Theta_4\neq0\},$

$\qquad U_{5}=\{L(\alpha )\in FLeib_{5}:\alpha _{3}=0,\Delta
_{4}\neq 0,\Theta_4=0\},$

$\qquad U_{6}=\{L(\alpha )\in FLeib_{5}:\alpha _{3}=0,\Delta
_{4}=0,\Theta_4\neq0\},$

$\qquad U_{7}=\{L(\alpha )\in FLeib_{5}:\alpha _{3}=0,\Delta
_{4}=0,\Theta_4=0\}.$\\

Now we will investigate the isomorphism problem for each of these
sets separately.

\textbf{Proposition 2.1.1} Two algebras $L(\alpha )$ and $L(\alpha
^{\prime })$ from $U_1$ are isomorphic if and only if

\[
\left( \frac{\alpha _{3}}{\Delta_4}\right)
^{2}\Theta_4=\left(\frac{\alpha^{\prime}_3}{\Delta_4'}\right)^{2}\Theta_4'
\]

This means that the expression $$\left( \frac{\alpha
_{3}}{\Delta_4}\right) ^{2}\Theta_4$$ can be taken as a parameter
$\lambda$ and then algebras from the set
$U_{1}$ can be parameterized as $L(1,0,\lambda ).$\\

\textbf{Proposition 2.1.2.}

a) All algebras from the set $U_2$ are isomorphic to $L(1,-2,0);$

 b) All algebras from the set $U_{3}$ are isomorphic to $L(1,-2,-2);$

 c) All algebras from the set $U_{4}$ are isomorphic to $L(0,1,0);$

 d) All algebras from the set $U_{5}$ are isomorphic to $L(0,1,1);$

 e) All algebras from the set $U_{6}$ are isomorphic to $L(0,0,1);$

 f) All algebras from the set $U_{7}$ are isomorphic to $L(0,0,0).$\\

\textbf{Theorem 2.1.3.} Let $L$ be a none Lie complex filiform
Leibniz algebra in $FLeib_{5}$.  Then it is isomorphic to one of
the following pairwise non-isomorphic Leibniz algebras:

$1)\ \ L(1,0,{\lambda}):$

$\qquad \qquad e_{0}e_{0}=e_{2},\ \ e_{i}e_{0}=e_{i+1,}\ \ 1\leq
i\leq3,\ \ e_{0}e_{1}=e_{3}+\lambda e_{4}, \ \ e_{1}e_{1}=e_{3},$

$\qquad \qquad e_{2}e_{1}=e_{4},\ \ \lambda \in
\mathbf{C}.$\\

$2) \ \ L(1,-2,0):$

$\qquad \qquad e_{0}e_{0}=e_{2},\ \ e_{i}e_{0}=e_{i+1,}\ \ 1\leq
i\leq3,\ \ e_{0}e_{1}=e_{3},\ \ e_{1}e_{1}=e_{3}-2e_{4},$

$\qquad \qquad e_{2}e_{1}=e_{4}.$\\

$3) \ \ L(1,-2,-2):$

$\qquad \qquad e_{0}e_{0}=e_{2},\ \ e_{i}e_{0}=e_{i+1,}\ \ 1\leq
i\leq3 ,\ \ e_{0}e_{1}=e_{3}-2e_{4}, $

$\qquad \qquad e_{1}e_{1}=e_{3}-2e_{4},\ \ e_{2}e_{1}=e_{4}
.$\\

$4) \ \ L(0,1,0):$

$\qquad \qquad e_{0}e_{0}=e_{2},\ \ e_{i}e_{0}=e_{i+1,},\ \ 1\leq
i\leq3 ,\ \ e_{1}e_{1}=e_{4}. $\\

$5) \ \ L(0,1,1):$

$\qquad \qquad e_{0}e_{0}=e_{2},\ \ e_{i}e_{0}=e_{i+1,},\ \ 1\leq
i\leq3 ,\ \ e_{0}e_{1}=e_{4},\ \ e_{1}e_{1}=e_{4}.$\\

$6) \ \ L(0,0,1):$

$\qquad \qquad e_{0}e_{0}=e_{2},\ \ e_{i}e_{0}=e_{i+1,},\ \ 1\leq
i\leq3 ,\ \ e_{0}e_{1}=e_{4}.$\\

$7) \ \ L(0,0,0):$

$\qquad \qquad e_{0}e_{0}=e_{2},\ \ e_{i}e_{0}=e_{i+1,}\ \ 1\leq i\leq3.$%
\\

The number of isomorphism classes $N_5=7.$ \\

\textbf{2.2 Dimension 6}

 Now we consider the six dimensional case. The set $FLeib_6$ can
be represented as a disjoint union of the subsets:

$\qquad FLeib_{6}=U_{1}\bigcup U_{2}\bigcup U_{3}\bigcup
U_{4}\bigcup U_{5}\bigcup U_{6}\bigcup U_{7}\bigcup U_{8}\bigcup
U_{9}\bigcup U_{10}\bigcup U_{11},$ where

$\qquad U_1=\{L(\alpha)\in FLeib_6:\alpha_3\neq 0, \Delta_4 \neq 0
\},$

$\qquad U_{2}=\{L(\alpha )\in FLeib_{6}:\alpha _{3}\neq
0,\Delta_4=0,\Delta_5 \neq 0,\Theta_5\neq0\},$

$\qquad U_{3}=\{L(\alpha )\in FLeib_{6}:\alpha _{3}=0,\Delta
_{4}\neq 0,\Delta _{5}\neq 0\},$

$\qquad U_{4}=\{L(\alpha )\in FLeib_{6}:\alpha _{3}\neq
0,\Delta_4=0,\Delta_5 \neq 0,\Theta_5=0\},$

$\qquad U_{5}=\{L(\alpha )\in FLeib_{6}:\alpha _{3}\neq
0,\Delta_4=0,\Delta_5=0,\Theta_5\neq0\},$

$\qquad U_{6}=\{L(\alpha )\in FLeib_{6}:\alpha _{3}=0,\Delta
_{4}\neq 0,\Delta _{5}=0,\Theta_5\neq 0\},$

$\qquad U_{7}=\{L(\alpha )\in FLeib_{6}:\alpha _{3}=0,\Delta
_{4}\neq 0,\Delta_{5}=0,\Theta_5 =0\},$

$\qquad U_{8}=\{L(\alpha )\in FLeib_{6}:\alpha _{3}=0,\Delta
_{4}=0,\Delta _{5}\neq 0,\Theta_5\neq0\},$

$\qquad U_{9}=\{L(\alpha )\in FLeib_{6}:\alpha _{3}=0,\Delta
_{4}=0,\Delta _{5}\neq 0,\Theta_5=0\},$

$\qquad U_{10}=\{L(\alpha )\in FLeib_{6}:\alpha _{3}=0,\Delta
_{4}=0,\Delta _{5}=0,\Theta_5 \neq 0\},$

$\qquad U_{11}=\{L(\alpha )\in FLeib_{6}:\alpha _{3}=0,\Delta
_{4}=0,\Delta_{5}=0,\Theta_5 =0\}.$\newline

\textbf{Proposition 2.2.1.} Two algebras $L(\alpha )$ and
$L(\alpha ^{\prime })$ from $U_1$ are isomorphic if and only if

\[
\frac{\alpha _{3}(\Delta_{5}+5\alpha
_{3}\Delta_4)}{\Delta_4^{2}}=\frac{\alpha
_{3}'(\Delta_{5}'+5\alpha _{3}'\Delta_4')}{\Delta_4'^{2}}
\]
\[
\frac{\alpha _{3}^{3}\Theta_{5}}{\Delta_4^{3}}=\frac{\alpha
_{3}'^{3}\Theta_{5}'}{\Delta_4'^{3}}
\]

Thus the following two expressions can be taken as parameters
$\lambda_{1},\lambda_{2}:$
$$\frac{\alpha _{3}(\Delta_{5}+5\alpha
_{3}\Delta_4)}{\Delta_4^{2}},$$
$$\frac{\alpha _{3}^{3}\Theta_{5}}{\Delta_4^{3}}$$ and algebras from $U_{1}$ can
be parameterized as $$L(1,0,{\lambda }_{1},{\lambda }_{2}).$$

\textbf{Proposition 2.2.2.}  Two algebras $L(\alpha )$ and
$L(\alpha ^{\prime })$ from $U_{2}$ are isomorphic if and only if

\[
\frac{\Delta_5^{3}}{\alpha _{3}^{3}\Theta_5^{2}}=\
\frac{\Delta_5'^{3}}{\alpha _{3}'^{3}\Theta_5'^{2}}.
\]

Thus in the set $U_{2}$ the expression
$$\frac{\Delta_5^{3}}{\alpha _{3}^{3}\Theta_5^{2}}$$ can be taken
as a parameter and algebras from $U_{2}$ can be parameterized as
$$L(1,-2,{\lambda },2{\lambda }-5).$$

\textbf{Proposition 2.2.3.}  Two algebras $L(\alpha )$ and
$L(\alpha ^{\prime })$ from $ U_{3} $ are isomorphic if and only
if

\[
\frac{\alpha _{4}^{3}\Theta_5}{\alpha _{5}^{3}}=\frac{\alpha
_{4}'^{3}\Theta_5'}{\alpha _{5}'^{3}}.
\]

So as a parameter $\lambda$ in $U_{3}$ we will take the expression
$$\frac{\alpha _{4}^{3}\Theta_5}{\alpha _{5}^{3}}$$ and write algebras from the set $U_{3}$ as
$$L(0,1,1,{\lambda }).$$

\textbf{Proposition 2.2.4.}

 a) All algebras from the set $U_{4}$ are isomorphic to the
algebra $L(1,-2,0,0);$

 b) All algebras from the set $U_5$ are isomorphic to the algebra $%
L(1,-2,5,0);$

 c) All algebras from the set $U_{5}$ are isomorphic to the
algebra $L(0,1,0,1);$

 d) All algebras from the set $U_{6}$ are isomorphic to the
algebra $L(0,1,0,0);$

 e) All algebras from the set $U_{7}$ are isomorphic to the
algebra $L(0,0,1,0);$

 f) All algebras from the set $U_{8}$ are isomorphic to the
algebra $L(0,0,1,1);$

 g) All algebras from the set $U_{9}$ are isomorphic to the
algebra $L(0,0,0,1);$

 h) All algebras from the set $U_{10}$ are isomorphic to
the algebra $L(0,0,0,0).$\\

\textbf{Theorem 2.2.5.} Let $L$ be a none Lie complex filiform
Leibniz algebra in $FLeib_{6}$. Then it is isomorphic to one of
the following pairwise non-isomorphic Leibniz algebras:

$1)\ \ L(1,0,{\lambda} _{1},{\lambda} _{2}):$

$\qquad \qquad e_{0}e_{0}=e_{2},\ \ e_{i}e_{0}=e_{i+1,}\ \ 1\leq
i\leq 4,\ \ e_{0}e_{1}=e_{3}+\lambda _{2}e_{5},$

$\qquad \qquad e_{1}e_{1}=e_{3}+{\lambda }_{1}e_{5},\ \
e_{2}e_{1}=e_{4},\ \ e_{3}e_{1}=e_{5},\ \ {\lambda }_{1},{\lambda
}_{2}\in \mathbf{C}.$\\

$2)\ \ L(1,-2,{\lambda} ,2{\lambda} -5):$

$\qquad \qquad e_{0}e_{0}=e_{2},\ \ e_{i}e_{0}=e_{i+1},\ \ 1\leq
i\leq 4,\ \ e_{0}e_{1}=e_{3}-2e_{4}+(2{\lambda }-5)e_{5},$

$\qquad \qquad e_{1}e_{1}=e_{3}-2e_{4}+{\lambda }e_{5},\ \
e_{2}e_{1}=e_{4}-2e_{5},\ \ e_{3}e_{1}=e_{5},\ \ {\lambda }\in \mathbf{C}.$%
\\

$3)\ \ L(0,1,1,{\lambda} ):$

$\qquad \qquad e_{0}e_{0}=e_{2},\ \ e_{i}e_{0}=e_{i+1},\ \ 1\leq
i\leq 4,\ \ e_{0}e_{1}=e_{4}+{\lambda }e_{5},$

$\qquad \qquad e_{1}e_{1}=e_{4}+e_{5},\ \ e_{2}e_{1}=e_{5},\ \
{\lambda }\in \mathbf{C}.$\\

$4)\ \ L(1,-2,0,0):$

$\qquad \qquad e_{0}e_{0}=e_{2},\ \ e_{i}e_{0}=e_{i+1},\ \ 1\leq
i\leq 4,\ \ e_{0}e_{1}=e_{3}-2e_{4},$

$\qquad \qquad e_{1}e_{1}=e_{3}-2e_{4}, \ \
e_{2}e_{1}=e_{4}-2e_{5},\ \ e_{3}e_{1}=e_{5}.$\\

$5)\ \ L(1,-2,5,0):$

$\qquad \qquad e_{0}e_{0}=e_{2},\ \ e_{i}e_{0}=e_{i+1},\ \ 1\leq
i\leq 4,\ \ e_{0}e_{1}=e_{3}-2e_{4},$

$\qquad \qquad e_{1}e_{1}=e_{3}-2e_{4}+5e_{5},\ \
e_{2}e_{1}=e_{4}-2e_{5},\ \ e_{3}e_{1}=e_{5}.$\\

$6)\ \ L(0,1,0,1):$

$\qquad \qquad e_{0}e_{0}=e_{2},\ \ e_{i}e_{0}=e_{i+1},\ \ 1\leq
i\leq 4,\ \ e_{0}e_{1}=e_{4}+e_{5},\ \ e_{1}e_{1}=e_{4},$

$\qquad \qquad e_{2}e_{1}=e_{5}.$\\

$7)\ \ L(0,1,0,0):$

$\qquad \qquad e_{0}e_{0}=e_{2},\ \ e_{i}e_{0}=e_{i+1},\ \ 1\leq
i\leq 4,\ \ e_{0}e_{1}=e_{4},\ \ e_{1}e_{1}=e_{4},\ \
e_{2}e_{1}=e_{5}.$\\

$8)\ \ L(0,0,1,0):$

$\qquad \qquad e_{0}e_{0}=e_{2},\ \ e_{i}e_{0}=e_{i+1},\ \ 1\leq
i\leq 4,\ \ e_{1}e_{1}=e_{5}.$\\

$9)\ \ L(0,0,1,1):$

$\qquad \qquad e_{0}e_{0}=e_{2},\ \ e_{i}e_{0}=e_{i+1},\ \ 1\leq
i\leq 4,\ \ e_{0}e_{1}=e_{5},\ \ e_{1}e_{1}=e_{5}.$\\

$10)\ \ L(0,0,0,1):$

$\qquad \qquad e_{0}e_{0}=e_{2},\ \ e_{i}e_{0}=e_{i+1},\ \ 1\leq
i\leq 4,\ \ e_{0}e_{1}=e_{5}.$\\

$11)\ \ L(0,0,0,0):$

$\qquad \qquad e_{0}e_{0}=e_{2},\ \ e_{i}e_{0}=e_{i+1}, \ \ 1\leq i\leq 4.$%
\\

The number of isomorphism classes $N_6=11.$ \\

\textbf{2.3 Dimension 7}

$FLeib_{7}=U_{1}\bigcup U_{2}\bigcup U_{3}\bigcup U_{4}\bigcup
U_{5}\bigcup U_{6}\bigcup U_{7}\bigcup U_{8}\bigcup U_{9}\bigcup
U_{10}\bigcup U_{11}\bigcup U_{12}\bigcup $

$\qquad \qquad \  U_{13}\bigcup U_{14}\bigcup U_{15}\bigcup
U_{16}\bigcup U_{17}, $

 where

$U_{1}=\{L(\alpha )\in FLeib_{7}:\alpha _{3}\neq 0,\Delta_4\neq
0\},$

$U_{2}=\{L(\alpha )\in FLeib_{7}:\alpha _{3}\neq 0,\Delta_4=0,
\Delta_5\neq 0,\Delta_{6}+6\alpha _{3}\Delta_5\neq 0\},$

$U_{3}=\{L(\alpha )\in FLeib_{7}:\alpha _{3}\neq 0,\Delta_4=0,
\Delta_5\neq 0, \Delta_{6}+6\alpha _{3}\Delta_5=0\},$

$U_{4}=\{L(\alpha )\in FLeib_{7}:\alpha _{3}\neq
0,\Delta_4=0,\Delta_{5}=0,\Delta_6 \neq 0,\Theta _{6}\neq 0\},$

$U_{5}=\{L(\alpha )\in FLeib_{7}:\alpha _{3}=0,\Delta_{4}\neq
0,\Delta_{5}\neq 0\},$

$U_{6}=\{L(\alpha )\in FLeib_{7}:\alpha _{3}=0,\Delta_{4}\neq
0,\Delta_{5}=0,\Delta_{6}+3\Delta_{4}^{2}\neq 0,\Theta_{6}\neq
0\},$

$ U_{7}=\{L(\alpha )\in FLeib_{7}:\alpha _{3}=0,\Delta
_{4}=0,\Delta_{5}\neq 0,\Delta_{6}\neq 0\},$

$U_{8}=\{L(\alpha )\in FLeib_{7}:\alpha _{3}\neq
0,\Delta_4=0,\Delta_5=0, \Delta_6\neq 0,\Theta _{6}=0\},$

$U_{9}=\{L(\alpha )\in FLeib_{7}:\alpha _{3}\neq
0,\Delta_4=0,\Delta_5=0, \Delta_6=0,\Theta_{6}\neq 0\},$

$ U_{10}=\{L(\alpha )\in FLeib_{7}:\alpha _{3}=0,\Delta_{4}\neq
0,\Delta_{5}=0,\Delta_{6}+3\Delta_{4}^{2}\neq 0,\Theta_{6}=0\},$

$U_{11}=\{L(\alpha )\in FLeib_{7}:\alpha _{3}=0,\Delta_{4}\neq
0,\Delta_{5}=0,\Delta_{6}+3\Delta_{4}^{2}=0,\Theta_{6}\neq 0\},$

 $U_{12}=\{L(\alpha )\in FLeib_{7}:\alpha _{3}=0,\Delta
_{4}=0,\Delta_{5}\neq 0,\Delta_{6}= 0,\Theta_6 \neq 0\},$

$ U_{13}=\{L(\alpha )\in FLeib_{7}:\alpha _{3}=0,\Delta
_{4}=0,\Delta _{5}\neq 0,\Delta_{6}= 0,\Theta_6 =0\},$

$ U_{14}=\{L(\alpha )\in FLeib_{7}:\alpha _{3}=0,\Delta
_{4}=0,\Delta_{5}=0,\Delta_{6}\neq 0,\Theta_{6}\neq 0\},$

$U_{15}=\{L(\alpha )\in FLeib_{7}:\alpha _{3}=0,\Delta
_{4}=0,\Delta_{5}=0,\Delta_{6}\neq 0,\Theta_{6}=0\},$

$U_{16}=\{L(\alpha )\in FLeib_{7}:\alpha _{3}=0,\Delta
_{4}=0,\Delta_{5}=0,\Delta_{6}=0,\Theta_6 \neq 0\},$

$ U_{17}=\{L(\alpha )\in FLeib_{7}:\alpha _{3}=0,\Delta
_{4}=0,\Delta _{5}=0,\Delta_{6}=0,\Theta_6 =0\}.$\\

\textbf{Proposition 2.3.1.} Two algebras $L(\alpha )$ and
$L(\alpha')$ from $U_1$ are isomorphic if and only if

\[
\frac{\alpha _{3}(\Delta_{5}+5\alpha
_{3}\Delta_4)}{\Delta_4^{2}}=\frac{\alpha
_{3}'(\Delta_{5}'+5\alpha _{3}'\Delta_4')}{\Delta_4'^{2}}
\]
\[
\frac{\alpha
_{3}(\alpha_{3}\Delta_{6}+6\alpha_{3}^{2}\Delta_{5}-3\Delta_{4}\Delta_{5}+9\alpha_{3}^{3}\Delta_{4}-12\alpha_{3}
\Delta_{4}^{2})}{\Delta_4^{3}}=\]
\[ \frac{\alpha
_{3}'(\alpha_{3}'\Delta_{6}'+6\alpha_{3}'^{2}\Delta_{5}'-3\Delta_{4}'\Delta_{5}'+9\alpha_{3}'^{3}\Delta_{4}'-12\alpha_{3}'
\Delta_{4}'^{2})}{\Delta_{4}'^{3}}
\]

\[
\frac{\alpha _{3}^4\Theta_6}{\Delta_{4}^4}=\frac{\alpha
_{3}'^4\Theta_6'}{\Delta_{4}'^4}.
\]

Thus, in this case algebras from the set $U_{1}$ can be
parameterized as $L(1,0,\lambda _{1},\lambda _{2},\lambda
_{3}).$\\

\textbf{Proposition 2.3.2} Two algebras $L(\alpha )$ and
$L(\alpha')$ from $U_{2}$ are isomorphic if and only if

\[
\frac{\Delta_5^{3}}{\alpha_3(\Delta_6+6\alpha_3\Delta_5)^2}=\frac{\Delta_5'^{3}}{\alpha_3'(\Delta_6'+6\alpha_3'\Delta_5')^2}
\]
\begin{eqnarray*}
\frac{\Delta_5^{4}\Theta_{6}}{\left(\Delta_6+6\alpha_3\Delta_5\right)
^{4}}
=\frac{\Delta_5'^{4}\Theta_{6}'}{\left(\Delta_6'+6\alpha_3'\Delta_5'\right)
^{4}}
\end{eqnarray*}

The expressions the above can be taken as parameters in $U_{2}$
and the set $U_{2}$ can be represented as $L(1,-2,\lambda _{1},
-5\lambda _{1}-14,\lambda _{2}).$\\

\textbf{Proposition 2.3.3.} Two algebras $L(\alpha )$ and
$L(\alpha')$ from $U_{3}$ are isomorphic if and only if

\[
\frac{\alpha _{3}^{4}\Theta_{6}^{2}}{\Delta_{5}^{4}}=\frac{\alpha
_{3}'^{4}\Theta_{6}'^{2}}{\Delta_{5}'^{4}}.
\]

The parameter $\lambda$ for algebras from the set $U_{3}$ is
$$\frac{\alpha _{3}^{4}\Theta_{6}}{\Delta_{5}^{4}}$$ and $U_{3}$ can be parameterized as
$L(1,-2,0,16,\lambda ).$\\

\textbf{Proposition 2.3.4.} Two algebras $L(\alpha )$ and
$L(\alpha')$ from $U_{4}$ are isomorphic if and only if

\[
\frac{\Delta_6^{4}}{\alpha
_{3}^{4}\Theta_{6}^{3}}=\frac{\Delta_6'^{4}}{\alpha
_{3}'^{4}\Theta_{6}'^{3}}.
\]

$U_{4}$ can be parameterized as $L(1,-2,5,\lambda ,2\lambda -14).$\\

\textbf{Proposition 2.3.5.} Two algebras $L(\alpha )$ and
$L(\alpha')$ from $U_{5}$ are isomorphic if and only if

\[
\frac{\Delta_{4}(\Delta_{6}+3\Delta_{4}^{2})}{\Delta
_{5}^{2}}=\frac{\Delta_{4}'(\Delta_{6}'+3\Delta_{4}'^{2})}{\Delta
_{5}'^{2}},
\]

\[
\left(\frac{\Delta
_{4}}{\Delta_{5}}\right)^4\Theta_{6}=\left(\frac{\Delta
_{4}'}{\Delta_{5}'}\right)^4\Theta_{6}'.
\]

The set $U_{5}$ can be parameterized as $L(0,1,1,\lambda _{1},\lambda _{2}).$\\

\textbf{Proposition 2.3.6.} Two algebras $L(\alpha )$ and
$L(\alpha')$ from $U_{6}$ are isomorphic if and only if

\[
\frac{(\Delta_{6}+3\Delta_{4}^{2})^2}{\Delta_{4}^{2}\Theta_{6}}=\frac{(\Delta_{6}'+3\Delta_{4}'^{2})^2}{\Delta_{4}'^{2}\Theta_{6}'}.
\]

$U_{6}$ can be represented as a parameterized family of algebras
$L(0,1,0,\lambda ,2\lambda -3).$\\

\textbf{Proposition 2.3.7.} Two algebras $L(\alpha )$ and
$L(\alpha')$ from $U_{7}$ are isomorphic if and only if

\[
\left(\frac{\Delta_{5}}{\Delta_{6}}\right)^4\Theta_6=\left(\frac{\Delta_{5}'}{\Delta_{6}'}\right)^4\Theta_6'.
\]

We will get one parametric family of algebras for the set $U_{7}:$
$L(0,0,1,1,\lambda ).$\\

\textbf{Proposition 2.3.8}

a) All algebras from the set $U_{8}$ are isomorphic to
$L(1,-2,5,0,0);$

b) All algebras from $U_{9}$ are isomorphic to $L(1,-2,5,14,0);$

 c) Algebras from $U_{10}$ are isomorphic to $L(0,1,0,0,0);$

 d) All algebras from $U_{11}$ are isomorphic to $L(0,1,0,-3,0);$

 e) All algebras from $U_{12}$ are isomorphic to $L(0,0,1,0,1);$

 f) All algebras from $U_{13}$ are isomorphic to $L(0,0,1,0,0);$

 g) All algebras from $U_{14}$ are isomorphic to $L(0,0,0,1,0);$

 h) Algebras from $U_{15}$ are isomorphic to $L(0,0,0,1,1);$

 i) Algebras from $U_{16}$ are isomorphic to $L(0,0,0,0,1);$

 j) Algebras from $U_{17}$ are isomorphic to $L(0,0,0,0,0).$\\

\textbf{Theorem 2.3.9.} Let $L$ be  a none Lie complex filiform
Leibniz algebra in $FLeib_{7}$. Then it is isomorphic to one of
the following pairwise non-isomorphic Leibniz algebras:

$1)\ \ L(1,0,{\lambda} _{1},{\lambda} _{2},{\lambda} _{3}):$

$\qquad \qquad e_{0}e_{0}=e_{2},\ \ e_{i}e_{0}=e_{i+1},\ \ 1\leq
i\leq 5,\ \ e_{0}e_{1}=e_{3}+{\lambda }_{1}e_{5}+{\lambda
}_{3}e_{6},$

$\qquad \qquad e_{1}e_{1}=e_{3}+{\lambda} _{1}e_{5}+{\lambda}
_{2}e_{6},\ \ e_{2}e_{1}=e_{4}+{\lambda} _{1}e_{6},\ \
e_{3}e_{1}=e_{5},\ \ e_{4}e_{1}=e_{6},$

$\qquad \qquad {\lambda }_{1},{\lambda }_{2},{\lambda }_{3}\in \mathbf{C}.$%
\\

$2)\ \ L(1,-2,{\lambda} _{1},-5{\lambda} _{1}-14,{\lambda} _{2}):$

$\qquad \qquad e_{0}e_{0}=e_{2},\ \ e_{i}e_{0}=e_{i+1},\ \ 1\leq
i\leq 5,\ \ e_{0}e_{1}=e_{3}-2e_{4}+{\lambda }_{1}e_{5}+{\lambda
}_{2}e_{6},$

$\qquad \qquad e_{1}e_{1}=e_{3}-2e_{4}+{\lambda}
_{1}e_{5}+(-5{\lambda} _{1}-14)e_{6},\ \
e_{2}e_{1}=e_{4}-2e_{5}+{\lambda} _{1}e_{6},$

$\qquad \qquad e_{3}e_{1}=e_{5}-2e_{6},\ \ e_{4}e_{1}=e_{6},\ \ {\lambda }%
_{1},{\lambda }_{2}\in \mathbf{C}.$\\

$3)\ \ L(1,-2,0,16,{\lambda} ):$

$\qquad \qquad e_{0}e_{0}=e_{2},\ \ e_{i}e_{0}=e_{i+1},\ \ 1\leq
i\leq 5,\ \ e_{0}e_{1}=e_{3}-2e_{4}+{\lambda }e_{6},$

$\qquad \qquad e_{1}e_{1}=e_{3}-2e_{4}+16e_{6},\ \
e_{2}e_{1}=e_{4}-2e_{5},\ \ e_{3}e_{1}=e_{5}-2e_{6},$

$\qquad \qquad e_{4}e_{1}=e_{6},\ \ {\lambda }\in
\mathbf{C}.$\\

$4)\ \ L(1,-2,5,{\lambda},2{\lambda} -14):$

$\qquad \qquad e_{0}e_{0}=e_{2},\ \ e_{i}e_{0}=e_{i+1},\ \ 1\leq
i\leq 5,\ \ e_{0}e_{1}=e_{3}-2e_{4}+5e_{5}+(2{\lambda }-14)e_{6},$

$\qquad \qquad e_{1}e_{1}=e_{3}-2e_{4}+5e_{5}+{\lambda} e_{6},\ \
e_{2}e_{1}=e_{4}-2e_{5}+5e_{6},\ \ e_{3}e_{1}=e_{5}-2e_{6},$

$\qquad \qquad e_{4}e_{1}=e_{6},\ \ {\lambda }\in
\mathbf{C}.$\\

$5)\ \ L(0,1,1,{\lambda} _{1},{\lambda} _{2}):$

$\qquad \qquad e_{0}e_{0}=e_{2},\ \ e_{i}e_{0}=e_{i+1},\ \ 1\leq
i\leq 5,\ \ e_{0}e_{1}=e_{4}+e_{5}+{\lambda }_{2}e_{6},$

$\qquad \qquad e_{1}e_{1}=e_{4}+e_{5}+{\lambda }_{1}e_{6},\ \
e_{2}e_{1}=e_{5}+e_{6},\ \ e_{3}e_{1}=e_{6},\ \ {\lambda }_{1},{\lambda }%
_{2}\in \mathbf{C}.$\\

$6)\ \ L(0,1,0,{\lambda} ,2{\lambda}-3): $

$\qquad \qquad e_{0}e_{0}=e_{2},\ \ e_{i}e_{0}=e_{i+1},\ \ 1\leq
i\leq 5,\ \ e_{0}e_{1}=e_{4}+(2{\lambda }-3)e_{6},$

$\qquad \qquad e_{1}e_{1}=e_{4}+{\lambda }e_{6},\ \
e_{2}e_{1}=e_{4},\ \ e_{3}e_{1}=e_{6},\ \ {\lambda }\in
\mathbf{C}.$\\

$7)\ \ L(0,0,1,1,{\lambda}):$

$\qquad \qquad e_{0}e_{0}=e_{2},\ \ e_{i}e_{0}=e_{i+1},\ \ 1\leq
i\leq 5,\ \ e_{0}e_{1}=e_{5}+{\lambda }e_{6},\ \
e_{1}e_{1}=e_{5}+e_{6},$

$\qquad \qquad e_{2}e_{1}=e_{6},\ \ {\lambda }\in \mathbf{C}.\
$\\

$8)\ \ L(1,-2,5,0,0):$

$\qquad \qquad e_{0}e_{0}=e_{2},\ \ e_{i}e_{0}=e_{i+1},\ \ 1\leq
i\leq 5,\ \ e_{0}e_{1}=e_{3}-2e_{4}+5e_{5},$

$\qquad \qquad e_{1}e_{1}=e_{3}-2e_{4}+5e_{5},\ \
e_{2}e_{1}=e_{4}-2e_{5}+5e_{6},\ \ e_{3}e_{1}=e_{5}-2e_{6},\ \
e_{4}e_{1}=e_{6}.$\\

$9)\ \ L(1,-2,5,14,0):$

$\qquad \qquad e_{0}e_{0}=e_{2},\ \ e_{i}e_{0}=e_{i+1},\ \ 1\leq
i\leq 5,\ \ e_{0}e_{1}=e_{3}-2e_{4}+5e_{5},$

$\qquad \qquad e_{1}e_{1}=e_{3}-2e_{4}+5e_{5}-14e_{6},\ \
e_{2}e_{1}=e_{4}-2e_{5}+5e_{6},\ \ e_{3}e_{1}=e_{5}-2e_{6},$

$\qquad \qquad e_{4}e_{1}=e_{6}.$\\

$10)\ \ L(0,1,0,0,0):$

$\qquad \qquad e_{0}e_{0}=e_{2},\ \ e_{i}e_{0}=e_{i+1},\ \ 1\leq
i\leq 5,\ \ e_{0}e_{1}=e_{4},\ \ e_{1}e_{1}=e_{4},\ \
e_{2}e_{1}=e_{5},$

$\qquad \qquad e_{3}e_{1}=e_{6}.$\\

$11)\ \ L(0,1,0,-3,0):$

$\qquad \qquad e_{0}e_{0}=e_{2},\ \ e_{i}e_{0}=e_{i+1},\ \ 1\leq
i\leq 5,\ \ e_{0}e_{1}=e_{4},\ \ e_{1}e_{1}=e_{4}-3e_{6},$

$\qquad \qquad e_{2}e_{1}=e_{5},\ \ e_{3}e_{1}=e_{6}.$\\

$12)\ \ L(0,0,1,0,1):$

$\qquad \qquad e_{0}e_{0}=e_{2},\ \ e_{i}e_{0}=e_{i+1,}\ \ 1\leq
i\leq 5,\ \ e_{0}e_{1}=e_{5}+e_{6},\ \ e_{2}e_{1}=e_{6}.$\\
$13)\ \ L(0,0,1,0,0):$

$\qquad \qquad e_{0}e_{0}=e_{2},\ \ e_{i}e_{0}=e_{i+1},\ \ 1\leq
i\leq 5,\ \ e_{0}e_{1}=e_{5},\ \ e_{2}e_{1}=e_{6}.$\\

$14)\ \ L(0,0,0,1,0):$

$\qquad \qquad e_{0}e_{0}=e_{2},\ \ e_{i}e_{0}=e_{i+1,}\ \ 1\leq
i\leq 5,\ \ e_{1}e_{1}=e_{6}.$\\

$15)\ \ L(0,0,0,1,1):$

$\qquad \qquad e_{0}e_{0}=e_{2},\ \ e_{i}e_{0}=e_{i+1},\ \ 1\leq
i\leq 5,\ \ e_{0}e_{1}=e_{6,}\ \ e_{1}e_{1}=e_{6}.$\\

$16)\ \ L(0,0,0,0,1):$

$\qquad \qquad e_{0}e_{0}=e_{2},\ \ e_{i}e_{0}=e_{i+1,}\ \ 1\leq
i\leq 5,\ \ e_{0}e_{1}=e_{6}.$\\

$17)\ \ L(0,0,0,0,0):$

$\qquad \qquad e_{0}e_{0}=e_{2},\ \ e_{i}e_{0}=e_{i+1},\ \ 1\leq i\leq 5.$%
\\

The number of isomorphism classes $N_7=17.$ \\

\textbf{2.4 Dimension 8}

 $FLeib_{8}=U_{1}\bigcup U_{2}\bigcup U_{3}\bigcup
U_{4}\bigcup U_{5}\bigcup U_{6}\bigcup U_{7}\bigcup U_{8}\bigcup
U_{9}\bigcup U_{10}\bigcup U_{11}\bigcup U_{12}\bigcup U_{13}
\bigcup$

$\qquad \qquad \ \ U_{14}\bigcup U_{15}\bigcup U_{16}\bigcup
U_{17}\bigcup U_{18}\bigcup U_{19}\bigcup U_{20}\bigcup \bigcup
U_{21}\bigcup U_{22}\bigcup U_{23}\bigcup U_{24}\bigcup U_{25},$

where

$U_{1}=\{L(\alpha )\in FLeib_{8}:\alpha _{3}\neq 0,\Delta _{4}\neq
0\},$

$U_{2}=\{L(\alpha )\in FLeib_{8}:\alpha _{3}\neq 0,\Delta
_{4}=0,\Delta_{5}\neq 0,\Delta_{6}+6\alpha _{3}\Delta_{5}\neq
0\},$

$U_{3}=\{L(\alpha )\in FLeib_{8}:\alpha _{3}\neq 0,\Delta
_{4}=0,\Delta_{5}\neq 0,\Delta_{6}+6\alpha _{3}\Delta_{5}\ =0,
\Theta_{7}\neq 0\},$

$U_{4}=\{L(\alpha )\in FLeib_{8}:\alpha _{3}\neq 0,\Delta
_{4}=0,\Delta_{5}\neq 0, \Delta_{6}+6\alpha _{3}\Delta_{5}\
=0,\Theta_{7}= 0\},$

$U_{5}=\{L(\alpha )\in FLeib_{8}:\alpha _{3}\neq 0,\Delta
_{4}=0,\Delta_{5}=0,\Delta_{6}\neq 0, \Delta_{7}+7\alpha
_{3}\Delta_{6}\neq 0\},$

$U_{6}=\{L(\alpha )\in FLeib_{8}:\alpha _{3}\neq
0,\Delta_{4}=0,\Delta_{5}=0, \Delta _{6}\neq 0, \Delta_{7}+7\alpha
_{3}\Delta_{6}=0\},$

$U_{7}=\{L(\alpha )\in FLeib_{8}:\alpha _{3}\neq 0,\Delta
_{4}=0,\Delta_{5}=0, \Delta_{6}=0, \Delta_{7}\neq 0\},$

$U_{8}=\{L(\alpha )\in FLeib_{8}:\alpha _{3}=0,\Delta_{4}\neq
0,\Delta_{5}\neq 0\},$

$ U_{9}=\{L(\alpha )\in FLeib_{8}:\alpha _{3}=0,\Delta_{4}\neq
0,\Delta_{5}=0,\Delta_{6}+3\Delta_{4}^{2}\neq 0,\Delta _{7}\neq
0\},$

$ U_{10}=\{L(\alpha )\in FLeib_{8}:\alpha _{3}=0,\Delta _{4}\neq
0,\Delta_{5}=0,\Delta_{6}+3\Delta_{4}^{2}\neq 0,\Theta _{7} \neq
0\},$

$ U_{11}=\{L(\alpha )\in FLeib_{8}:\alpha _{3}=0,\Delta_{4}\neq
0,\Delta_{5}=0,\Delta_{6}+3\Delta_{4}^{2}=0,\Delta _{7}\neq
0,\Theta_7\neq 0\},$

$ U_{12}=\{L(\alpha )\in FLeib_{8}:\alpha _{3}=0,\Delta
_{4}=0,\Delta_{5}\neq 0,\Delta_{6}\neq 0\},$

$U_{13}=\{L(\alpha )\in FLeib_{8}:\alpha _{3}=0,\Delta
_{4}=0,\Delta_{5}\neq 0,\Delta_{6}= 0,\Delta_{7}\neq 0,\Theta
_{7}\neq 0\},$

$U_{14}=\{L(\alpha )\in FLeib_{8}:\alpha _{3}=0,\Delta
_{4}=0,\Delta_{5}=0,\Delta_{6}\neq 0,\Delta_{7}\neq 0\},$

$ U_{15}=\{L(\alpha )\in FLeib_{8}:\alpha _{3}=0,\Delta_{4}\neq
0,\Delta_{5}=0,\Delta_{6}+3\Delta_{4}^{2}\neq 0,\Delta
_{7}=0,\Theta_7=0\}, $

$U_{16}=\{L(\alpha )\in FLeib_{8}:\alpha _{3}=0,\Delta_{4}\neq
0,\Delta_{5}=0,\Delta_{6}+3\Delta_{4}^{2}=0,\Delta_{7}\neq
0,\Theta_7= 0\},$

$U_{17}=\{L(\alpha )\in FLeib_{8}:\alpha _{3}=0,\Delta_{4}\neq
0,\Delta_{5}=0,\Delta_{6}+3\Delta_{4}^{2}=0,\Delta_{7}=0,\Theta_7
\neq 0\},$

$U_{18}=\{L(\alpha )\in FLeib_{8}:\alpha _{3}=0,\Delta
_{4}=0,\Delta_{5}\neq 0,\Delta_{6}= 0,\Delta_{7}\neq
0,\Theta_{7}=0\},$

$U_{19}=\{L(\alpha )\in FLeib_{8}:\alpha _{3}=0,\Delta
_{4}=0,\Delta_{5}\neq 0,\Delta_{6}= 0,\Delta_{7}=0,\Theta_7\neq
0\},$

 $U_{20}=\{L(\alpha )\in FLeib_{8}:\alpha _{3}=0,\Delta
_{4}=0,\Delta_{5}=0,\Delta_{6}\neq 0,\Delta_{7}=0,\Theta_7 \neq
0\},$

$U_{21}=\{L(\alpha )\in FLeib_{8}:\alpha _{3}=0,\Delta
_{4}=0,\Delta_{5}=0,\Delta_{6}\neq 0,\Delta_{7}=0,\Theta_7 =0\},$

$U_{22}=\{L(\alpha )\in FLeib_{8}:\alpha _{3}=0,\Delta
_{4}=0,\Delta_{5}=0,\Delta_{6}=0,\Delta_{7}\neq 0,\Theta_7\neq
0\},$

$U_{23}=\{L(\alpha )\in FLeib_{8}:\alpha _{3}=0,\Delta
_{4}=0,\Delta _{5}=0,\Delta_{6}=0,\Delta_{7}\neq 0,\Theta_7=0\},$

 $U_{24}=\{L(\alpha )\in FLeib_{8}:\alpha
_{3}=0,\Delta_{4}=0,\Delta_{5}=0,\Delta_{6}=0,\Delta_{7}=0,\Theta_7
\neq 0\},$

$U_{25}=\{L(\alpha )\in FLeib_{8}:\alpha _{3}=0,\Delta
_{4}=0,\Delta_{5}=0,\Delta_{6}=0,\Delta_{7}=0,\Theta_7=0\}.$\\

\textbf{Proposition 2.4.1.} Two algebras $L(\alpha )$ and
$L(\alpha')$ from $U_1$ are isomorphic if and only if

\[
\frac{\alpha _{3}(\Delta_{5}+5\alpha
_{3}\Delta_4)}{\Delta_4^{2}}=\frac{\alpha
_{3}'(\Delta_{5}'+5\alpha _{3}'\Delta_4')}{\Delta_4'^{2}}
\]
\[
\frac{\alpha
_{3}(\alpha_{3}\Delta_{6}+6\alpha_{3}^{2}\Delta_{5}-3\Delta_{4}\Delta_{5}+9\alpha_{3}^{3}\Delta_{4}-12\alpha_{3}
\Delta_{4}^{2})}{\Delta_4^{3}}= \]
\[ \frac{\alpha
_{3}'(\alpha_{3}'\Delta_{6}'+6\alpha_{3}'^{2}\Delta_{5}'-3\Delta_{4}'\Delta_{5}'+9\alpha_{3}'^{3}\Delta_{4}'-12\alpha_{3}'
\Delta_{4}'^{2})}{\Delta_{4}'^{3}}
\]

\[
\frac{\alpha _{3}^{3}\Delta_{7}+28\alpha _{3}^{4}\Delta
_{4}^{2}+7\alpha _{3}^{6}\Delta_4+14\alpha _{3}^{5}\Delta
_{5}+7\alpha _{3}^{4}\Delta_6+7\alpha
_{3}^{3}\Delta_{4}\Delta_{5}}{\Delta_{4}^{4}}= \] \qquad \qquad
\qquad \[ \frac{\alpha _{3}'^{3}\Delta_{7}'+28\alpha
_{3}'^{4}\Delta _{4}'^{2}+7\alpha _{3}'^{6}\Delta_4'+14\alpha
_{3}'^{5}\Delta _{5}'+7\alpha _{3}'^{4}\Delta_6'+7\alpha
_{3}'^{3}\Delta_{4}'\Delta_{5}'}{\Delta_{4}'^{4}}
\]

\[
\frac{\alpha _{3}^{5}\Theta_{7}}{\Delta_{4}^{5}} =\frac{\alpha
_{3}'^{5}\Theta_{7}'}{\Delta_{4}'^{5}}
\]

$U_{1}$ is parameterized as
$L(1,0,\lambda _{1},\lambda _{2},\lambda _{3},\lambda _{4}).$\\

\textbf{Proposition 2.4.2.} Two algebras $L(\alpha )$ and
$L(\alpha')$ from $U_{2}$ are isomorphic if and only if

\[
\frac{\Delta_{5}^{3}}{\alpha _{3}(\Delta_{6}+6\alpha
_{3}\Delta_{5})^2}=\frac{\Delta_{5}'^{3}}{\alpha
_{3}'(\Delta_{6}'+6\alpha _{3}'\Delta_{5}')^2}
\]
\[
\frac{\Delta_{5}^{4}(\Delta_{7}+7\alpha _{3}\Delta_{6}+14\alpha
_{3}^2\Delta_{5})}{\alpha _{3}(\Delta_{6}+6\alpha _{3}\Delta
_{5})^{4}}=\frac{\Delta_{5}'^{4}(\Delta_{7}'+7\alpha
_{3}'\Delta_{6}'+14\alpha _{3}'^2\Delta_{5}')}{\alpha
_{3}'(\Delta_{6}'+6\alpha _{3}'\Delta _{5}')^{4}}
\]
\[
\frac{\Delta_{5}^{5}\Theta_{7}}{(\Delta_{6}+6\alpha _{3}\Delta
_{5})^{5}} =\frac{\Delta_{5}'^{5}\Theta_{7}'}{(\Delta_{6}'+6\alpha
_{3}'\Delta _{5}')^{5}}
\]

The set $U_{2}$ can be parameterized as $L(1,-2,\lambda
_{1},-5\lambda _{1}-14,\lambda _{2},\lambda _{3}).
$\\

\textbf{Proposition 2.4.3.} Two algebras $L(\alpha )$ and
$L(\alpha')$ from $U_{3}$ are isomorphic if and only if

\[
\left( \frac{\Delta_{5}}{\alpha
_{3}}\right)^{5}\frac{1}{\Theta_{7}^{2}}=\left(\frac{\Delta_{5}'}{\alpha
_{3}}'\right)^{5}\frac{1}{\Theta_{7}'^{2}}
\]
\[
\frac{\Delta_{5}^{8}(\Delta_{7}-28\alpha _{3}^{2}\Delta
_{5})}{\alpha _{3}^{9}\Theta_{7}^{4}}
=\frac{\Delta_{5}'^{8}(\Delta_{7}'-28\alpha _{3}'^{2}\Delta
_{5}')}{\alpha _{3}'^{9}\Theta_{7}'^{4}}
\]

The algebras from the set $U_{3}$ can be parameterized as
$L(1,-2,\lambda _{1},-6\lambda _{1}-14,\lambda _{2},\lambda
_{1}^{2}).$\\

\textbf{Proposition 2.4.4.} Two algebras $L(\alpha )$ and
$L(\alpha')$ from $U_{4}$ are isomorphic if and only if

\[
\frac{\alpha _{3}(\Delta_{7}-28\alpha _{3}^{2}\Delta_{5})}{\Delta
_{5}^{2}}=\frac{\alpha _{3}'(\Delta_{7}'-28\alpha
_{3}'^{2}\Delta_{5}')}{\Delta _{5}'^{2}}.
\]

The set $U_4$ can be parameterized as
$L(1,-2,0,16,\lambda ,\lambda ).$\\

 \textbf{Proposition 2.4.5.} Two algebras $L(\alpha )$ and $L(\alpha')$ from $U_{5}$
 are isomorphic if and only if

\[
\frac{\Delta_{6}^{4}}{\alpha _{3}(\Delta_{7}+7\alpha _{3}\Delta
_{6})^{3}}\\
=\frac{\Delta_{6}'^{4}}{\alpha _{3}'(\Delta_{7}'+7\alpha
_{3}'\Delta _{6}')^{3}}.
\]
\[
\left(\frac{\Delta_{6}}{\Delta_{7}+7\alpha _{3}\Delta
_{6}}\right)^5\Theta_7
=\left(\frac{\Delta_{6}'}{\Delta_{7}'+7\alpha _{3}'\Delta
_{6}'}\right)^5\Theta_7'
\]

$U_{5}$ is parameterized as
$L(1,-2,5,\lambda _{1},-6\lambda _{1}+42,\lambda _{2}).$\\

\textbf{Proposition 2.4.6.} Two algebras $L(\alpha )$ and
$L(\alpha')$ from $U_{6}$ are isomorphic if and only if

\[
\left(\frac{\alpha
_{3}}{\Delta_{6}}\right)^{5}\Theta_{7}^{3}=\left(\frac{\alpha
_{3}'}{\Delta_{6}'}\right)^{5}\Theta_{7}'^{3}
\]

The set $U_{6}$ can be parameterized as
$L(1,-2,5,\lambda ,-7\lambda +42,\lambda ^{2}).$\\

\textbf{Proposition 2.4.7.} Two algebras $L(\alpha )$ and
$L(\alpha ^{\prime })$ from $U_{7}$ are isomorphic if and only if

\[
\left(\frac{\Delta_{7}}{\alpha
_{3}}\right)^5\frac{1}{\Theta_{7}^{4}}=\left(\frac{\Delta_{7}'}{\alpha
_{3}'}\right)^5\frac{1}{\Theta_{7}'^{4}}
\]%

$U_{7}$ is parameterized as
$L(1,-2,5,-14,\lambda,2({\lambda} +21)).$\\

\textbf{Proposition 2.4.8.} Two algebras $L(\alpha )$ and
$L(\alpha')$ from $U_{8} $ are isomorphic if and only if
\[
\frac{\Delta_{4}(\alpha _{6}+3\alpha _{4}^{2})}{\Delta
_{5}^{2}}=\frac{\Delta_{4}'(\alpha _{6}'+3\alpha
_{4}'^{2})}{\Delta _{5}'^{2}}
\]

\[
\frac{\Delta_{4}^{2}(\Delta_{7}+7\Delta_{4}\Delta_{5})}{\Delta
_{5}^{3}}=\frac{\Delta_{4}'^{2}(\Delta
_{7}'+7\Delta_{4}'\Delta_{5}')}{\Delta _{5}'^{3}}.
\]
\[
\left(\frac{\Delta_{4}}{\Delta_{5}}\right)^{5}\Theta_{7}=\left(\frac{\Delta_{4}'}{\Delta_{5}'}\right)^{5}\Theta_{7}'.
\]

$U_8$ can be parameterized
as $L(0,1,1,\lambda _{1},\lambda _{2},\lambda _{3}).$\\

\textbf{Proposition 2.4.9.} Two algebras $L(\alpha )$ and
$L(\alpha ^{\prime })$ from $U_{9} $ are isomorphic if and only if

\[
\frac{(\Delta_{6}+3\Delta_{4}^{2})^{3}}{\Delta_{4}\Delta
_{7}^{2}}=\frac{(\Delta_{6}'+3\Delta_{4}'^{2})^{3}}{\Delta_{4}'\Delta
_{7}'^{2}}
\]%
\[
\left(\frac{\Delta_{6}+3\Delta
_{4}^{2}}{\Delta_{7}}\right)^{5}\Theta_{7}=\left(\frac{\Delta_{6}'+3\Delta
_{4}'^{2}}{\Delta_{7}'}\right)^{5}\Theta_{7}'.
\]

Thus the algebras from the set $U_{9}$ can be parameterized as
$L(0,1,0,\lambda _{1},\lambda _{1}+3,\lambda _{2}).$\\

\textbf{Proposition 2.4.10.} Two algebras $L(\alpha )$ and
$L(\alpha ^{\prime })$ from $U_{10}$ are isomorphic if and only if

\[
\left(\frac{\Delta_{6}+3\Delta_{4}^{2}}{\Delta_{4}}\right)^{5}\Theta_7^2
=\left(\frac{\Delta_{6}'+3\Delta_{4}'^{2}}{\Delta_{4}'}\right)^{5}\Theta_7'^{2}.
\]

$U_{10}$ can be parameterized as
$L(0,1,0,\lambda ,0,\lambda ^{2}).$\\

\textbf{Proposition 2.4.11.} Two algebras $L(\alpha )$ and $L(\alpha ^{\prime })$ from $%
U_{11}$ are isomorphic if and only if

\[
\left(\frac{\Delta_{4}}{\Delta_{7}}\right)^{5}\Theta
_{7}^{3}=\left(\frac{\Delta_{4}'}{\Delta_{7}'}\right)^{5}\Theta
_{7}'^{3}.
\]

$L(0,1,0,-3,\lambda ,\lambda ^{2}+\lambda )$ are representatives
of $U_{11}.$\\

\textbf{Proposition 2.4.12.} Two algebras $L(\alpha )$ and
$L(\alpha ^{\prime })$ from $U_{12}$ are isomorphic if and only if

\[
\frac{\Delta_{5}\Delta_{7}}{\Delta
_{6}^{2}}=\frac{\Delta_{5}'\Delta_{7}'}{\Delta _{6}'^{2}}
\]
\[
\left(\frac{\Delta_{5}}{\Delta_{6}}\right)^{5}\Theta_{7}=\left(\frac{\Delta_{5}'}{\Delta_{6}'}\right)^{5}\Theta_{7}'.
\]

Thus, the algebras from the set $U_{12}$ can be parameterized as
$L(0,0,1,1,\lambda _{1},\lambda _{2}).$\\

\textbf{Proposition 2.4.13.} Two algebras $L(\alpha )$ and
$L(\alpha ^{\prime })$ from $U_{13}$ are isomorphic if and only if

\[
\left(\frac{\Delta_{7}}{\Delta_{5}}\right)^{5}\frac{1}{\Theta_{7}^{2}}=\left(\frac{\Delta_{7}'}{\Delta_{5}'}\right)^{5}\frac{1}{\Theta_{7}'^{2}}.
\]

Thus, the algebras from the set $U_{13}$ can be parameterized as
$L(0,0,1,0,\lambda ,\lambda ^{2}+\lambda ).$\\

\textbf{Proposition 2.4.14.} Two algebras $L(\alpha )$ and
$L(\alpha ^{\prime })$ from $U_{14}$ are isomorphic if and only if

\[
\left(\frac{\Delta_{6}}{\Delta_{7}}\right)^{5}\Theta_{7}=\left(\frac{\Delta_{6}'}{\Delta_{7}'}\right)^{5}\Theta_{7}'.
\]

$L(0,0,0,1,1,\lambda )$ is a parametrization of $U_{14}.$\\

\textbf{Proposition 2.4.15}

a) All algebras from the set $U_{15}$ are isomorphic to
$L(0,1,0,0,0,0);$

b) All algebras from the set $U_{16}$ are isomorphic to
$L(0,1,0,-3,1,1);$

c) All algebras from the set $U_{17}$ are isomorphic to
$L(0,1,0,-3,0,1);$

d)All algebras from the set $U_{18}$ are isomorphic to
$L(0,0,1,0,1,1);$

 e) All algebras from the set $U_{19}$ are isomorphic to $L(0,0,1,0,0,1);$

 f) All algebras from the set $U_{20}$ are isomorphic to $L(0,0,0,1,0,1);$

g) All algebras from the set $U_{21}$ are isomorphic to
$L(0,0,0,1,0,0);$

h) All algebras from the set $U_{22}$ are isomorphic to
$L(0,0,0,0,1,0);$

i) All algebras from the set $U_{23}$ are isomorphic to
$L(0,0,0,0,1,1);$

j) All algebras from the set $U_{24}$ are isomorphic to
$L(0,0,0,0,0,1);$

 k) All algebras from the set $U_{25}$ are isomorphic to
 $L(0,0,0,0,0,0).$\\

\textbf{Theorem 2.4.16.} Let $L$ be a none Lie complex filiform
Leibniz algebra in $FLeib_{8}$. Then it is isomorphic to one of
the following pairwise non-isomorphic Leibniz algebras:

$1)\ \ L(1,0,{\lambda} _{1},{\lambda} _{2},{\lambda}
_{3},{\lambda} _{4}):$

$\qquad \qquad e_{0}e_{0}=e_{2},\ \ e_{i}e_{0}=e_{i+1,}\ \ 1\leq i
\leq 6,\ \ e_{0}e_{1}=e_{3}+{\lambda} _{1}e_{5}+{\lambda}
_{2}e_{6}+{\lambda} _{4}e_{7},$

$\qquad \qquad e_{1}e_{1}=e_{3}+{\lambda} _{1}e_{5}+{\lambda} _{2}e_{6}+{%
\lambda} _{3}e_{7},\ \ e_{2}e_{1}=e_{4}+{\lambda} _{1}e_{6}+{\lambda}%
_{3}e_{7},$

$\qquad \qquad e_{3}e_{1}=e_{5}+{\lambda }_{1}e_{7},\ \
e_{4}e_{1}=e_{6},\ \
e_{5}e_{1}=e_{7},\ \ {\lambda }_{1},{\lambda }_{2},{\lambda }_{3},{\lambda }%
_{4}\in \mathbf{C}.$\\

$2)\ \ L(1,-2,{\lambda} _{1},-(5{\lambda} _{1}+14),{\lambda}
_{2},{\lambda} _{3}):$

$\qquad \qquad e_{0}e_{0}=e_{2},\ \ e_{i}e_{0}=e_{i+1,}\ \ 1\leq i
\leq 6,\
\ e_{0}e_{1}=e_{3}-2e_{4}+{\lambda} _{1}e_{5}+{\lambda}_{2}e_{6}+{\lambda}%
_{3}e_{7},$

$\qquad \qquad e_{1}e_{1}=e_{3}-2e_{4}+{\lambda}
_{1}e_{5}-(5{\lambda} _{1}+14)e_{6}+{\lambda} _{2}e_{7},$

$\qquad \qquad e_{2}e_{1}=e_{4}-2e_{5}+{\lambda}
_{1}e_{6}-(5{\lambda} _{1}+14)e_{7},\ \
e_{3}e_{1}=e_{5}-2e_{6}+{\lambda}_{1}e_{7},$

$\qquad \qquad e_{4}e_{1}=e_{6}-2e_{7},\ \ e_{5}e_{1}=e_{7},\ \ {\lambda }%
_{1},{\lambda }_{2},{\lambda }_{3}\in \mathbf{C}.$\\

$3)\ \ L(1,-2,{\lambda} _{1},-(6{\lambda} _{1}+14),{\lambda}
_{2},{\lambda} _{1}^{2}):$

$\qquad \qquad e_{0}e_{0}=e_{2},\ \ e_{i}e_{0}=e_{i+1,}\ \ 1\leq i
\leq 6,\ \ e_{0}e_{1}=e_{3}-2e_{4}+{\lambda} _{1}e_{5}+{\lambda}
_{2}e_{6}+{\lambda} _{1}^{2}e_{7},$

$\qquad \qquad e_{1}e_{1}=e_{3}-2e_{4}+{\lambda}
_{1}e_{5}-(6{\lambda} _{1}+14)e_{6}+{\lambda} _{2}e_{7},$

$\qquad \qquad e_{2}e_{1}=e_{4}-2e_{5}+{\lambda}
_{1}e_{6}-(6{\lambda} _{1}+14)e_{7},\ \
e_{3}e_{1}=e_{5}-2e_{6}+{\lambda} _{1}e_{7},$

$\qquad \qquad e_{4}e_{1}=e_{6}-2e_{7},\ \ e_{5}e_{1}=e_{7},\ \ {\lambda }%
_{1},{\lambda }_{2}\in \mathbf{C}.$\\

$4)\ \ L(1,-2,0,16,{\lambda} ,{\lambda} ):$

$\qquad \qquad e_{0}e_{0}=e_{2},\ \ e_{i}e_{0}=e_{i+1},\ \ 1\leq i
\leq 6,\ \ e_{0}e_{1}=e_{3}-2e_{4}+16e_{6}+{\lambda} e_{7},$

$\qquad \qquad e_{1}e_{1}=e_{3}-2e_{4}+16e_{6}+{\lambda} e_{7},\ \
e_{2}e_{1}=e_{4}-2e_{5}+16e_{7},$

$\qquad \qquad e_{3}e_{1}=e_{5}-2e_{6},\ \
e_{4}e_{1}=e_{6}-2e_{7},\ \ e_{5}e_{1}=e_{7},\ \ {\lambda }\in
\mathbf{C}.$\\

$5)\ \ L(1,-2,5,{\lambda}_{1},-6({\lambda} _{1}-7),{\lambda}
_{2}):$

$\qquad \qquad e_{0}e_{0}=e_{2},\ \ e_{i}e_{0}=e_{i+1},\ \ 1\leq i
\leq 6,\ \ e_{0}e_{1}=e_{3}-2e_{4}+5e_{5}+{\lambda}
_{1}e_{6}+{\lambda} _{2}e_{7},$

$\qquad \qquad e_{1}e_{1}=e_{3}-2e_{4}+5e_{5}+{\lambda} _{1}e_{6}-6({\lambda}%
_{1}-7)e_{7},$

$\qquad \qquad e_{2}e_{1}=e_{4}-2e_{5}+5e_{6}+{\lambda}
_{1}e_{7},\ \ e_{3}e_{1}=e_{5}-2e_{6}+5e_{7},$

$\qquad \qquad e_{4}e_{1}=e_{6}-2e_{7},\ \ e_{5}e_{1}=e_{7},\ \ {\lambda }%
_{1},{\lambda }_{2}\in \mathbf{C}.$\\

$6)\ \ L(1,-2,5,{\lambda} ,-7({\lambda}-6),{\lambda} ^{2}):$

$\qquad \qquad e_{0}e_{0}=e_{2},\ \ e_{i}e_{0}=e_{i+1},\ \ 1\leq i
\leq 6,\ \ e_{0}e_{1}=e_{3}-2e_{4}+5e_{5}+{\lambda}
e_{6}+{\lambda}^{2}e_{7},$

$\qquad \qquad e_{1}e_{1}=e_{3}-2e_{4}+5e_{5}+{\lambda}
e_{6}-7({\lambda}-6)e_{7},$

$\qquad \qquad e_{2}e_{1}=e_{4}-2e_{5}+5e_{6}+{\lambda} e_{7},\ \
e_{3}e_{1}=e_{5}-2e_{6}+5e_{7},$

$\qquad \qquad e_{4}e_{1}=e_{6}-2e_{7},\ \ e_{5}e_{1}=e_{7},\ \ {\lambda }%
\in \mathbf{C}.$\\

$7)\ \ L(1,-2,5,-14,{\lambda},2({\lambda} +21)):$

$\qquad \qquad e_{0}e_{0}=e_{2},\ \ e_{i}e_{0}=e_{i+1},\ \ 1\leq i
\leq 6,$

$\qquad \qquad  e_{0}e_{1}=e_{3}-2e_{4}+5e_{5}-14e_{6}+2({\lambda}
+21)e_{7}, $

$\qquad \qquad  e_{1}e_{1}=e_{3}-2e_{4}+5e_{5}-14e_{6}+{\lambda}
e_{7},\ \ e_{2}e_{1}=e_{4}-2e_{5}+5e_{6}-14e_{7},$

$\qquad \qquad e_{3}e_{1}=e_{5}-2e_{6}+5e_{7},\ \
e_{4}e_{1}=e_{6}-2e_{7},\ \ e_{5}e_{1}=e_{7},\ \ {\lambda}
\in \mathbf{C}.$\\

$8)\ \ L(0,1,1,{\lambda}_{1},{\lambda} _{2},{\lambda} _{3}):$

$\qquad \qquad e_{0}e_{0}=e_{2},\ \ e_{i}e_{0}=e_{i+1},\ \ 1\leq i
\leq 6,\ \ e_{0}e_{1}=e_{4}+e_{5}+{\lambda} _{1}e_{6}+{\lambda}
_{3}e_{7},$

$\qquad \qquad e_{1}e_{1}=e_{4}+e_{5}+{\lambda}
_{1}e_{6}+{\lambda} _{2}e_{7},\ \
e_{2}e_{1}=e_{5}+e_{6}+{\lambda}_{1}e_{7},$

$\qquad \qquad e_{3}e_{1}=e_{6}+e_{7},\ \ e_{4}e_{1}=e_{7},\ \ {\lambda }%
_{1},{\lambda }_{2},{\lambda }_{3}\in \mathbf{C}.$\\

$9)\ \ L(0,1,0,{\lambda}_{1},{\lambda}_{1}+3,{\lambda} _{2}):$

$\qquad \qquad e_{0}e_{0}=e_{2},\ \ e_{i1}e_{0}=e_{i+1,}\ \ 1\leq
i \leq 6,\ \ e_{0}e_{1}=e_{4}+{\lambda} _{1}e_{6}+{\lambda}
_{2}e_{7},$

$\qquad \qquad e_{1}e_{1}=e_{4}+{\lambda} _{1}e_{6}+({\lambda}
_{1}+3)e_{7},\ \ e_{2}e_{1}=e_{5}+{\lambda} _{1}e_{7},\ \
e_{3}e_{1}=e_{6},$

$\qquad \qquad e_{4}e_{1}=e_{7},\ \ {\lambda }_{1},{\lambda }_{2}\in \mathbf{%
C}.$\\

$10)\ \ L(0,1,0,{\lambda} ,0,{\lambda} ^{2}):$

$\qquad \qquad e_{0}e_{0}=e_{2},\ \ e_{i}e_{0}=e_{i+1,}\ \ 1\leq i
\leq 6,\ \ e_{0}e_{1}=e_{4}+{\lambda} e_{6}+{\lambda} ^{2}e_{7},$

$\qquad \qquad e_{1}e_{1}=e_{4}+{\lambda }e_{6},\ \ e_{2}e_{1}=e_{5}+{%
\lambda }e_{7},\ \ e_{3}e_{1}=e_{6},\ \ e_{4}e_{1}=e_{7},\ \
{\lambda }\in \mathbf{C}.$\\

$11)\ \ L(0,1,0,-3,{\lambda} ,{\lambda} ^{2}+{\lambda} ):$

$\qquad \qquad e_{0}e_{0}=e_{2},\ \ e_{i}e_{0}=e_{i+1,}\ \ 1\leq i
\leq6,\ \ e_{0}e_{1}=e_{4}-3e_{6}+({\lambda} ^{2}+{\lambda}
)e_{7},$

$\qquad \qquad e_{1}e_{1}=e_{4}-3e_{6}+{\lambda} e_{7},\ \
e_{2}e_{1}=e_{5}-3e_{7},\ \ e_{3}e_{1}=e_{6},\ \
e_{4}e_{1}=e_{7},$

$\qquad \qquad {\lambda }\in \mathbf{C}.$\\

$12)\ \ L(0,0,1,1,{\lambda}_{1},{\lambda} _{2}):$

$\qquad \qquad e_{0}e_{0}=e_{2},\ \ e_{i}e_{0}=e_{i+1},\ \ 1\leq i
\leq6, \ \ e_{0}e_{1}=e_{5}+e_{6}+{\lambda} _{2}e_{7},$

$\qquad \qquad e_{1}e_{1}=e_{5}+e_{6}+{\lambda }_{1}e_{7},\ \
e_{2}e_{1}=e_{6}+e_{7},\ \ e_{3}e_{1}=e_{7},\ \ {\lambda }_{1},{\lambda }%
_{2}\in \mathbf{C}.$\\

$13)\ \ L(0,0,1,0,{\lambda} ,{\lambda} ^{2}+{\lambda}):$

$\qquad \qquad e_{0}e_{0}=e_{2},\ \ e_{i}e_{0}=e_{i+1},\ \ 1\leq i
\leq6, \ \ e_{0}e_{1}=e_{5}+({\lambda}^{2}+{\lambda} )e_{7},$

$\qquad \qquad e_{1}e_{1}=e_{5}+{\lambda }e_{7},\ \
e_{2}e_{1}=e_{6},e_{3}e_{1}=e_{7},\ \ {\lambda }\in
\mathbf{C}.$\\

$14)\ \ L(0,0,0,1,1,{\lambda} ):$

$\qquad \qquad e_{0}e_{0}=e_{2},\ \ e_{i}e_{0}=e_{i+1},\ \ 1\leq i
\leq6,\ \ e_{0}e_{1}=e_{6}+{\lambda} e_{7},\ \
e_{1}e_{1}=e_{6}+e_{7},$

$\qquad \qquad e_{2}e_{1}=e_{7},\ \ {\lambda }\in
\mathbf{C}.$\\

$15)\ \ L(0,1,0,0,0,0):$

$\qquad \qquad e_{0}e_{0}=e_{2},\ \ e_{i}e_{0}=e_{i+1,}\ \ 1\leq i
\leq6,\ \ e_{0}e_{1}=e_{4},\ \ e_{1}e_{1}=e_{4},$

$\qquad \qquad e_{2}e_{1}=e_{5},\ \ e_{3}e_{1}=e_{6},\ \ e_{4}e_{1}=e_{7}.$%
\\

$16)\ \ L(0,1,0,-3,1,1):$

$\qquad \qquad e_{0}e_{0}=e_{2},\ \ e_{i}e_{0}=e_{i+1,}\ \ 1\leq i
\leq6,\ \ e_{0}e_{1}=e_{4}-3e_{6}+e_{7},$

$\qquad \qquad e_{1}e_{1}=e_{4}-3e_{6}+e_{7},\ \
e_{2}e_{1}=e_{5}-3e_{7},\ \ e_{3}e_{1}=e_{6},\ \
e_{4}e_{1}=e_{7}.$\\

$17)\ \ L(0,1,0,-3,0,1):$

$\qquad \qquad e_{0}e_{0}=e_{2},\ \ e_{i}e_{0}=e_{i+1},\ \ 1\leq i
\leq6, \ \ e_{0}e_{1}=e_{4}-3e_{6}+e_{7},$

$\qquad \qquad e_{1}e_{1}=e_{4}-3e_{6},\ \
e_{2}e_{1}=e_{5}-3e_{7},e_{3}e_{1}=e_{6},\ \
e_{4}e_{1}=e_{7}.$\\

$18)\ \ L(0,0,1,0,1,1):$

$\qquad \qquad e_{0}e_{0}=e_{2},\ \ e_{i}e_{0}=e_{i+1},\ \ 1\leq i
\leq6, \ \ e_{0}e_{1}=e_{5}+e_{7},$

$\qquad \qquad e_{1}e_{1}=e_{5}+e_{7},\ \ e_{2}e_{1}=e_{6},\ \
e_{3}e_{1}=e_{7}.$\\

$19)\ \ L(0,0,1,0,0,1):$

$\qquad \qquad e_{0}e_{0}=e_{2},\ \ e_{i}e_{0}=e_{i+1},\ \ 1\leq i
\leq6, \ \ e_{0}e_{1}=e_{5}+e_{7},\ \ e_{1}e_{1}=e_{5},$

$\qquad \qquad e_{2}e_{1}=e_{6},\ \ e_{3}e_{1}=e_{7}.$\\

$20)\ \ L(0,0,0,1,0,1):$

$\qquad \qquad e_{0}e_{0}=e_{2},\ \ e_{i}e_{0}=e_{i+1,}\ \ 1\leq i
\leq6,\ \ e_{0}e_{1}=e_{6}+e_{7},\ \ e_{1}e_{1}=e_{6},$

$\qquad \qquad e_{2}e_{1}=e_{7}.$\\

$21)\ \ L(0,0,0,1,0,0):$

$\qquad \qquad e_{0}e_{0}=e_{2},\ \ e_{i}e_{0}=e_{i+1}, \ \ 1\leq
i \leq6 ,\ \ e_{0}e_{1}=e_{6},\ \ e_{1}e_{1}=e_{6},\ \
e_{2}e_{1}=e_{7} .$\\

$22)\ \ L(0,0,0,0,1,0):$

$\qquad \qquad e_{0}e_{0}=e_{2},\ \ e_{i}e_{0}=e_{i+1},\ \ 1\leq i
\leq6,\ \ e_{1}e_{1}=e_{7}.$\\

$23)\ \ L(0,0,0,0,1,1):$

$\qquad \qquad e_{0}e_{0}=e_{2},\ \ e_{i}e_{0}=e_{i+1},\ \ 1\leq i
\leq6, \ \ e_{0}e_{1}=e_{7},\ \ e_{1}e_{1}=e_{7}.$\\

$24)\ \ L(0,0,0,0,0,1):$

$\qquad \qquad e_{0}e_{0}=e_{2},\ \ e_{i}e_{0}=e_{i+1},\ \ 1\leq i
\leq6, \ \ e_{0}e_{1}=e_{7}.$\\

$25)\ \ L(0,0,0,0,0,0):$

$\qquad \qquad e_{0}e_{0}=e_{2},\ \ e_{i}e_{0}=e_{i+1},\ \ 1\leq i
\leq6. $\\

The number of isomorphism classes $N_8=25.$ \\

\textbf{Conjecture.} The number of isomorphism classes $N_n$ of
$n$-dimensional none Lie complex filiform Leibniz algebras in
$FLeib_{n}$ can be found by the formula:
$$N_n=n^2-7n+17.$$

Note that the validity of the above formula is confirmed in
dimension $9$ as well.

We would like to thank U.D.Bekbaev and B.A.Omirov for their
helpful discussions.

\end{document}